\documentclass[journal]{IEEEtran}

\usepackage[T1]{fontenc}
\usepackage{graphicx}
\usepackage{epstopdf} 
\usepackage{amsmath,amsfonts,amssymb,amsthm,amstext,amscd,amsxtra,amsopn,stmaryrd,mathtools,mathrsfs}
\usepackage{cite}
\usepackage[]{subfigure}
\usepackage{multirow}
\usepackage{esint}
\usepackage{empheq}
\usepackage{bm}
\usepackage{algorithm} 
\usepackage{algpseudocode}
\usepackage{leftidx}
\usepackage{makecell}
\usepackage{stackrel}
\usepackage{multicol}

\usepackage{color,soul}
\soulregister\cite7
\soulregister\ref7
\soulregister\pageref7
\usepackage{slashbox}
\providecommand{\keywords}[1]{\textbf{\textit{Index terms---}} #1}

\usepackage{xcolor}

\usepackage{caption}

\title{A Fast Volume Integral Equation Solver \\ with Linear Basis Functions for the \\ Accurate Computation of EM  Fields in MRI}
\author{Ioannis P. Georgakis, Ilias I. Giannakopoulos, Mikhail S. Litsarev, \\ and Athanasios G. Polimeridis, ~\IEEEmembership{Senior Member,~IEEE}
	\thanks{This work was supported by grants from the Skoltech-MIT Next Generation Program.}
	\thanks{Ioannis P. Georgakis and Ilias I. Giannakopoulos are with the Skoltech Center for Computational and Data-Intensive Science and Engineering, Skolkovo Institute of Science and Technology, 121205 Moscow, Russia (e-mails: ioannis.georgakis@skolkovotech.ru, ilias.giannakopoulos@skolkovotech.ru).}
	\thanks{Mikhail S. Litsarev and Athanasios G. Polimeridis are with Q Bio, Redwood City, CA 94063, USA (emails: mchf.lms@gmail.com, thanospol88@gmail.com).}}

\begin{document}
	
	\onecolumn
	
	\noindent{\huge IEEE Copyright Notice}
	
	\vspace{0.5cm}
	\noindent{© 2020 IEEE.  Personal use of this material is permitted.  Permission from IEEE must be obtained for all other uses, in any current or future media, including reprinting/republishing this material for advertising or promotional purposes, creating new collective works, for resale or redistribution to servers or lists, or reuse of any copyrighted component of this work in other works.}
	
	\vspace{1.5cm}
	\noindent{ \large Accepted to be published in: 2020 IEEE Transactions on Antennas and Propagation.}
	
	\vspace{2.5cm}
	\noindent {Cite as:}
	
	\vspace{0.2cm}
	\noindent{I. P. Georgakis, I. I. Giannakopoulos, M. S. Litsarev and A. G. Polimeridis, "A Fast Volume Integral Equation Solver with Linear Basis Functions for the Accurate Computation of EM Fields in MRI," in IEEE Transactions on Antennas and Propagation, doi: 10.1109/TAP.2020.3044685.}
	\twocolumn

	\maketitle
	
	\begin{abstract}
		A stable volume integral equation (VIE) solver based on polarization/magnetization currents is presented, for the accurate and efficient computation of the electromagnetic scattering from highly inhomogeneous and high contrast objects. We employ the Galerkin Method of Moments to discretize the formulation with discontinuous piecewise linear basis functions on uniform voxelized grids, allowing for the acceleration of the associated matrix-vector products in an iterative solver, with the help of FFT. Numerical results illustrate the superior accuracy and more stable convergence properties of the proposed framework, when compared against standard low order (piecewise constant) discretization schemes and a more conventional VIE formulation based on electric flux densities. Finally, the developed solver is applied to analyze complex geometries, including realistic human body models, typically used in modeling the interactions between electromagnetic waves and biological tissue.
	\end{abstract}
	
	
	\section{INTRODUCTION}
	
	\IEEEPARstart{N}{umerical} modeling of electromagnetic (EM) scattering from extremely inhomogeneous objects and objects with high contrast is of great interest. Over the last decades, a plethora of numerical algorithms and computational methods has been developed for applications in wireless communications, microwave instrumentation, and remote sensing. However, more challenging problems arise, such as modeling the interactions between EM waves and biological tissue, which is of great relevance for the determination of the deposited radio-frequency (RF) power inside the human body and the associated safety considerations at high-field magnetic resonance imaging (MRI). Specifically, at high and ultra-high magnetic field strengths, the operating frequency of the associated RF coils also increases. Hence, the wavelength becomes comparable to or even smaller than the effective dimension of the human body, potentially increasing the specific absorption rate (SAR) and deteriorating the magnetic field homogeneity and thus the image quality. Therefore, the accurate estimation of the local distribution of the EM fields, generated by an external RF source, as well as the absorbed power, described by SAR, calls for the implementation of novel and more efficient electrodynamic solvers, to achieve a reliable numerical solution. Developing such precise simulation tools poses a great challenge since from a computational EM perspective the human body is highly inhomogeneous and presents great geometrical complexity.
	
	Regarding the numerical techniques for computing the EM fields in inhomogeneous dielectrics, such as the human body, there is a vast amount of literature and a rich investigative activity. On one hand, there are methods based on the partial differential form of Maxwell equations, such as the finite differences (FD)  and finite element (FE) method, which are known for their generality and versatility in many application areas. On the other hand, there are methods that start from Maxwell equations in their integro-differential form, such as the integral equation (IE) methods, which offer great flexibility in exploiting certain problem properties by customization of the method. IE methods, despite being complicated and computationally expensive, have proven to be the method of choice for modeling inhomogeneous objects by dividing the arbitrary scatterer into simple volumetric tessellations, with piecewise homogeneous properties. Furthermore, a volume integral equation (VIE)-based solver can be effortlessly embedded into a surface integral equation (SIE) solver to calculate the EM field distributions and RF transmit and receive coils interactions in the presence of realistic human body models (RHBMs) \cite{villena2016fast}. However, the development of an efficient and stable numerical EM simulation software based on integral equations is far from being a trivial task. Towards that direction, there has been a recent contribution \cite{polimeridis2014stable}, where a stable current-based VIE solver with piecewise constant (PWC) basis functions has been developed and incorporated into \cite{MARIE} for the fast EM analysis of MR coils.  
	
	In the MR modeling case, there is an internal limitation with respect to the spatial resolution and thus the grid refinements in a numerical solver. This naturally leads to the exploitation of $p$-refinement techniques for achieving more reliable solutions in a VIE-based solver, for coarse resolutions, too. This is particularly useful since, in many applications within the MR field, the models come in fixed resolutions without available refinements and it is important to compute fast, yet reliably, the EM fields distributions for personalized MRI applications such as the patient-specific SAR calculation and monitoring \cite{milshteyn2020individualized}. In addition, in costly inverse problems, that are computationally feasible for coarse resolutions due to the memory limitations of graphical  processing  units  (GPUs), $p$-refinement can provide more accurate solutions \cite{serralles2019noninvasive,serralles2017volumetric}. In this work, we modify the aforementioned current-based VIE solver \cite{polimeridis2014stable} and equip it with higher-order basis functions, that yield more accurate results for modeling highly inhomogeneous and high contrast objects, without the need of refining the mesh. Namely, discontinuous piecewise linear (PWL) \cite{tsai1986procedure,markkanen2014discretization} basis and testing functions are utilized and defined with support a single voxel of a uniform grid that allows for the fast calculation of the associated matrix-vector products with the help of FFT, exploiting the block Toeplitz structure of the system's matrix, since the arising integral kernels are translationally invariant. Plenty of research studies exist in the literature regarding this problem \cite{rubinacci2006broadband,lu2003fast,ozdemir2006nonconformal,ozdemir2007nonconformal,catedra1989numerical,gan1995discrete,sun2009novel,li2006applying,chew2001fast,jin1996computation,costabel2010volume,zhang2017discontinuous,borup1984fast,markkanen2012analysis,markkanen2012discretization,markkanen2016volume,yla2014surface,antenor1999scattering,schaubert1984tetrahedral,shen1989discrete,sancer2006volume,zwanborn1991weak,zwamborn1992three}, including an effective parallelization of the adaptive integral method \cite{wei20112,geyik2012fdtd,wei2013moreI,wei2013moreII} tailoring VIE for large-scale bioelectromagnetic applications and EM fields computations in anatomical human body models \cite{massey2016austinman}. Furthermore, both PWC and PWL methods constitute applications of the general $k$-space method \cite{bojarski1972k,bojarski1982k} for a specific choice of basis functions on a uniform cubic grid and the accuracy and convergence properties of the proposed method is expected to strongly depend on the value of the dielectric permittivity \cite{markkanen2014discretization}, similarly to the PWC case \cite{polimeridis2014stable}. The aforementioned currents must not necessarily satisfy any continuity conditions between neighboring elements, voxels in our case, allowing for the use of testing and basis functions that span $L^{2}$, when employing Method of Moments (MoM). This way, the finite-energy conditions are respected and the spectral properties of the operators are preserved, as it has been shown in recent studies \cite{van2007gaps,van2008well}. Moreover, the resulting formulation calls for the numerical integration of singular volume-volume Galerkin inner products, but as it has been shown in previous work \cite{georgakis2017reduction,polimeridis2013robust,bleszynski2013reduction,knockaert2011analytic}, there exist readily available formulas that reduce both the dimensionality and singularity order of the kernels, allowing for the fast and precise numerical evaluation of the singular integrals by means of well-established sophisticated cubatures \cite{polimeridis2008direct,polimeridis2010complete,polimeridis2011fast,polimeridis2013evaluation,polimeridis2013directfn,DEMCEM,DIRECTFN}, originally developed for SIE formulations. Recent contributions have also expanded the analysis for arbitrary quadrilateral patches \cite{tambova2017fully,tambova2018generalization}. 
	
	Finally, through numerical experiments, we illustrate that the resulting discretized formulation is well-conditioned and has more stable convergence properties than the discretized version with PWC approximations and than a more standard VIE formulation based on electric flux densities \cite{zwamborn1992three}. Specifically, the number of iterations of the iterative solver remains practically the same as with the case of PWC basis functions and is much smaller than that of the flux-based (DVIE) solver, for highly inhomogeneous scatterers.  Another favorable feature of the proposed scheme is its superior accuracy. To demonstrate that, we present a comparative analysis between the PWC and the PWL basis solver for a homogeneous sphere, comparing both solvers with analytical results obtained with Mie series \cite{mie1908beitrage}, and for RHBMs from the Virtual Family Population \cite{christ2009virtual}. Even when numerically treating the case of dielectric shimming \cite{van2017efficient}, where high-permittivity dielectric pads are placed in the vicinity of a RHBM, the suggested work performs effectively yielding EM field distributions comparable with the finite difference time domain (FDTD) method. FDTD is one of the most commonly used methods within the MR field, which, however, entails a large number of unknowns and long computational times and demands the application of absorbing boundary conditions (ABCs) that further expand the domain. IE methods, including our proposed method, come to the remedy of the abovementioned limitations since they automatically satisfy radiation conditions, thus there is no need for ABCs, are a highly customizable tool, and allow for higher-order approximations. Furthermore, IE methods, unlike FDTD, can be applied in the frequency domain which is suitable for MRI simulations that are single-frequency problems. However, PWL solver, compared against PWC, naturally, comes with an increase in the computational cost and in the required memory footprint. Thankfully, the arising Green function tensors in FFT-based VIEs preserve low multilinear rank properties \cite{polimeridis2014compression,giannakopoulos20183D}, thus,  we  can  overcome  this  limiting  factor and  accelerate  the  solution  of  the  linear  system  of  equations by combining tensor decompositions and parallel programming techniques of GPUs \cite{giannakopoulos2019memory}.

	The remainder of the paper is organized as follows. In Section II, we set up the current-based VIE formulation and describe its superior properties in comparison with other VIE formulations. In Section III, we introduce the novel discretization scheme with PWL basis functions, formulate the linear system by means of MoM, and represent the system in a tensor format for a better comprehension. In Section IV, we describe the procedure to accelerate the matrix-vector product and develop a fast FFT-based solver. Finally, in Section V, we validate the proposed solver by comparison with the analytical solution for spheres, study its convergence for dielectric cubes, and demonstrate its accuracy and convergence properties when modeling the EM scattering from RHBMs.
	
	\begin{table}[ht]
		\captionsetup{font=scriptsize}
		\caption{NOTATION} \label{table_notation} \centering
		\begin{tabular}{c| l }
			
			\hline\hline\\[-0.4em]
			Notation                              				  & Description                                                                \\[0.3em] \hline \\[-0.3em]
			$i$									  				  & imaginary unit															   \\[0.3em]
			$a$									  				  & scalar in $\mathbb{C}$												       \\[0.3em]
			$\bm{a}$                              				  & vector in $\mathbb{C}^3$,  $\bm{a} =(a_x,a_y,a_z)$                         \\[0.3em]
			$\mathbf{a}$                          				  & $1$-D array, vector in $\mathbb{C}^N$                                      \\[0.3em]
			$\mathbf{A}$                          				  & $2$-D array, matrix in $\mathbb{C}^{N_1\times N_2}$                        \\[0.3em]
			$\underline{\mathbf{A}}$			  				  & $m$-D array, tensor in $\mathbb{C}^{N_1\times N_2...\times N_m}$           \\[0.3em]
			$\mathcal{A}/\bm{\mathcal{A}}$		 			      & continuous operator acting in $\mathbb{C}^3/\mathbb{C}^6$                  \\[0.3em]
			$\odot{}$                             				  & element-wise multiplication                                                \\[0.3em]
			\hline\hline
		\end{tabular}
	\end{table}
	
	\section{VOLUME INTEGRAL EQUATION FORMULATION}
	We consider the scattering of time-harmonic EM waves by a penetrable object, occupying the bounded domain $\Omega$ in 3-D Euclidean space, $\mathbb{R}^{3}$. The working angular frequency is $\omega\in\mathbb{R}^{+}$ and the electrical properties are defined as
	\begin{equation}\label{eq1}
		\begin{split}
			\epsilon &= \epsilon_{0}, \quad \quad \quad \, \mu=\mu_{0} \quad \quad \quad \ \text{in} \quad \mathbb{R}^{3}\backslash \Omega, \\
			\epsilon &= \epsilon_{r}(\bm{r})\epsilon_{0}, \quad \mu=\mu_{r}(\bm{r})\mu_{0} \quad \text{in} \quad \Omega.
		\end{split}
	\end{equation}
	Here, the vacuum (or free-space) permittivity $\epsilon_{0}$ and permeability $\mu_{0}$ are real positive values, and the relative permittivities $\epsilon_{r}(\bm{r})$ and $\mu_{r}(\bm{r})$ are assumed complex:
	\begin{equation}\label{eq2}
		\begin{split}
			\epsilon_{r}(\bm{r}) &= \epsilon_{r}^{'}(\bm{r})-i\epsilon_{r}^{''}(\bm{r}),\\
			\mu_{r}(\bm{r}) &= \mu_{r}^{'}(\bm{r})-i\mu_{r}^{''}(\bm{r}),
		\end{split}
	\end{equation}
	with $\epsilon_{r}^{'}$, $\mu_{r}^{'} \in(0,\infty)$ and $\epsilon_{r}^{''}$, $\mu_{r}^{''} \in[0,\infty)$, assuming a time factor $\exp(i\omega t)$. At this point, we also define the free space wave number, $k_{0}=\omega\sqrt{\mu_{0}\epsilon_{0}}$.
	
	The total time harmonic fields $(\bm{e},\bm{h})$ in the presence of an isotropic and inhomogeneous object can be decomposed into incident $(\bm{e}_{{\rm inc}},\bm{h}_{{\rm inc}})$ and scattered $(\bm{e}_{{\rm sca}},\bm{h}_{{\rm sca}})$ fields:
	\begin{equation}\label{eq4}
		\left(\begin{array}{c}\bm{e}\\\bm{h}\end{array}\right)=\left(\begin{array}{c}\bm{e}_{{\rm inc}}\\ \bm{h}_{{\rm inc}}\end{array}\right)+\left(\begin{array}{c}\bm{e}_{{\rm sca}}\\\bm{h}_{{\rm sca}}\end{array}\right),
	\end{equation}
	where the incident fields are the fields generated by the impressed currents or the sources in the absence of the scatterer and the scattered fields are given by the induced currents due to the presence of the scatterer, as in {\cite{markkanen2012discretization,polimeridis2014stable,sun2009novel}}:
	\begin{equation}\label{eq13}
		\left(\begin{array}{c}\bm{e}_{{\rm sca}}\\\bm{h}_{{\rm sca}}\end{array}\right)=\left(\begin{array}{cc} {\frac{1}{c_{e}}(\mathcal{N}-\mathcal{I})} & -\mathcal{K}\\\mathcal{K} & {\frac{1}{c_{m}}(\mathcal{N}-\mathcal{I})}\end{array}\right)\left(\begin{array}{c}\bm{j}\\\bm{m}\end{array}\right),
	\end{equation}
	where the equivalent polarization and magnetization currents are defined as
	\begin{equation}\label{eq6}
		\begin{split}
			\bm{j}(\bm{r}) & \triangleq c_{e}\chi_{e}(\bm{r})\bm{e}(\bm{r}),\\
			\bm{m}(\bm{r}) & \triangleq c_{m}\chi_{m}(\bm{r})\bm{h}(\bm{r}),
		\end{split}
	\end{equation}
	with $c_{e},c_{m}\triangleq i\omega\epsilon_{0},i\omega\mu_{0}$ and the electric and magnetic susceptibilities are respectively given by
	\begin{equation}\label{eq7}
		\begin{split}
			\chi_{e}&\triangleq\epsilon_{r}-1,\\
			\chi_{m}&\triangleq\mu_{r}-1.
		\end{split}
	\end{equation}
	More explicitly, the associated continuous integro-differential operators are:
	\begin{subequations}\label{eq8}
		\begin{gather}
			\mathcal{N}\bm{u}\triangleq\nabla\times\nabla\times\mathcal{S}(\bm{u};\Omega)(\bm{r}),\\
			\mathcal{K}\bm{u} \triangleq\nabla\times\mathcal{S}(\bm{u};\Omega)(\bm{r}),
		\end{gather}
	\end{subequations}
	where
	\begin{equation}\label{eq9}
		\mathcal{S}(\bm{u};\Omega)(\bm{r})\triangleq\intop_{\Omega}g(\bm{r}-\bm{r}')\bm{u}(\bm{r}')d\bm{r}'
	\end{equation}
	and $g$ is the free-space scalar Green function:
	\begin{equation}\label{eq10}
		g(\bm{r})=\frac{e^{-ik_{0}|\bm{r}|}}{4\pi|\bm{r}|}.
	\end{equation}
	
	Finally, we derive the current-based volume integral equation for polarization and magnetization currents by combining (\ref{eq4}), (\ref{eq13}), and (\ref{eq6}), as in \cite{polimeridis2015computation}:
	\begin{equation}\label{eq14}
		\bm{\mathcal{A}}\left(\begin{array}{c}\bm{j}\\\bm{m}\end{array}\right)={\bm{\mathcal{C}}\bm{\mathcal{M}}_{\chi}}\left(\begin{array}{c}
			\bm{e}_{{\rm inc}}\\\bm{h}_{{\rm inc}}\end{array}\right),
	\end{equation}
	where
	\begin{subequations}\label{eq15}
		\begin{gather}
			\bm{\mathcal{A}}=\left(\begin{array}{cc}
				\mathcal{M}_{\epsilon_{r}}-\mathcal{M}_{\chi_{e}}\mathcal{N} & c_{e}\mathcal{M}_{\chi_{e}}\mathcal{K}\\
				-c_{m}\mathcal{M}_{\chi_{m}}\mathcal{K} & \mathcal{M}_{\mu_{r}}-\mathcal{M}_{\chi_{m}}\mathcal{N}
			\end{array}\right),\\
			\bm{\mathcal{C}}=\left(\begin{array}{cc}
				c_{e}\mathcal{I} & 0\\
				0 & c_{m}\mathcal{I}
			\end{array}\right),\\
			\bm{\mathcal{M}}_{\chi}=\left(\begin{array}{cc}
				\mathcal{M}_{\chi_{e}} & 0\\
				0 & \mathcal{M}_{\chi_{m}}
			\end{array}\right).
		\end{gather}
	\end{subequations}
	$\mathcal{M}_{\phi}$ is a multiplication operator that multiplies with the function $\phi$ and $\mathcal{I}$ is the identity dyadic tensor. The focal point of this work revolves around the interaction between EM waves and biological tissue, which can be considered purely dielectric, and (\ref{eq14}) reduces to the following:
	\begin{equation}\label{eq16}
		(\mathcal{M}_{\epsilon_{r}}-\mathcal{M}_{\chi_{e}}\mathcal{N})\bm{j}=c_{e}\mathcal{M}_{\chi_{e}}\bm{e}_{{\rm inc}},
	\end{equation}
	and the electric and magnetic field in the scatterer can be calculated as
	\begin{subequations}\label{tot_fields}
		\begin{gather}
			\bm{e}_{{\rm tot}} = \bm{e}_{{\rm inc}} + \frac{1}{c_{e}}(\mathcal{N}-\mathcal{I})\bm{j}, \\
			\bm{h}_{{\rm tot}} = \bm{h}_{{\rm inc}} + \mathcal{K}\bm{j}.
		\end{gather}
	\end{subequations}

	\section{DISCRETIZATION SCHEME}
	
	As demonstrated in recent studies, the current-based formulation has superior spectral properties when compared with the other alternatives, i.e., flux- and field-based formulations \cite{van2007gaps,van2008well,markkanen2012discretization,polimeridis2014stable}. In order to achieve this superior performance and guarantee convergence in the norm of the solution, when employing a Galerkin projection method, the function space of the basis and testing functions should be carefully chosen \cite{markkanen2016numerical,hsiao1997mathematical}. Specifically, the testing functions should span the $L^{2}$ dual of the range space of the associated operator, where the mapping properties of the current-based formulation, as expressed in (\ref{eq14}), read as:
	\begin{equation}
		\left[  {L}^2(\mathbb{R}^{3}) \right] ^3 \to \left[ {L}^2(\mathbb{R}^{3}) \right] ^3.
	\end{equation}
	Furthermore, as we describe in the next section, at the current-based formulation, the basis functions should not enforce any continuity conditions across the material interfaces and the discretization elements, so the natural function space of choice is that of square integrable functions $L^{2}$ \cite{van2007gaps,van2008well,markkanen2012discretization,polimeridis2014stable}. All these prerequisites are respected in a Galerkin projection method and are in accordance with previous work for a current-based VIE with PWC basis functions \cite{polimeridis2014stable}. However, the nature of the application at hand often requires higher accuracy, especially by means of $p$-refinement, i.e., using higher-order approximation of the VIE's unknown.
	
	\subsection{Piecewise Linear Basis Functions on Voxels}
	We begin our analysis by defining the computational domain and the associated discretization grid. We use a rectangular parallelepiped that encloses the inhomogeneous scatterer but does not have to extend any farther than that, since the radiation conditions at infinity are satisfied by the Green function. By denoting the dimensions of this rectangular box as $L_{x}$, $L_{y}$ and $L_{z}$, we can now discretize it with $N_{V}=N_{x}\times N_{y}\times N_{z}$
	number of voxels, which in the general case may be non-uniform and have side length $\Delta x=L_{x}/N_{x}$, $\Delta y=L_{y}/N_{y}$, and $\Delta z=L_{z}/N_{z}$. Having defined the grid, we expand the unknown equivalent polarization and magnetization currents in terms of suitable vector basis functions. These currents are unlikely to satisfy any continuity conditions because of the material discontinuities and as a result the function space of choice in a Galerkin discretization scheme should be $[{L}^{2}(\mathbb{R}^{3})]^{3}$. Without loss of generality, we proceed in the following with the expansion only for the polarization currents, since for the magnetization currents the basis functions would be identical, 
	\begin{equation}\label{eq18}
		\bm{j}(\bm{r})=j_{x}(\bm{r})\hat{\bm{x}}+j_{y}(\bm{r})\hat{\bm{y}}+j_{z}(\bm{r})\hat{\bm{z}},
	\end{equation}
	where each component of the current can be expanded in some discrete set of appropriate basis functions. In this work, we utilize the PWL basis functions, hence the current approximation reads
	\begin{equation}\label{eq19}
		j_{q}(\bm{r})\thickapprox\sum_{\bm{n}=(1,1,1)}^{(N_{x},N_{y},N_{z})}\sum_{l=1}^{4}N_{\bm{n}}^{l}(\bm{r})a_{\bm{n}}^{ql},
	\end{equation}
	where $q\in\{x,y,z\}$ indicates the components of the current, $\bm{n}=n_{x}\hat{\bm{x}}+n_{y}\hat{\bm{y}}+n_{z}\hat{\bm{z}}$ is a compound index denoting the centers of the voxels in the grid, $N_{\bm{n}}^{l}(\bm{r})$ are the $4$ basis functions per current component per voxel, and $a_{\bm{n}}^{ql}$ are the unknown coefficients. In detail, the scalar basis functions are defined as
	\begin{subequations}\label{eq20}
		\begin{gather}
			N_{\bm{n}}^{1}(\bm{r})=P_{\bm{n}}(\bm{r}),\\
			N_{\bm{n}}^{2}(\bm{r})=\frac{x-x_{n}}{\Delta x}P_{\bm{n}}(\bm{r}),\\
			N_{\bm{n}}^{3}(\bm{r})=\frac{y-y_{n}}{\Delta y}P_{\bm{n}}(\bm{r}),\\
			N_{\bm{n}}^{4}(\bm{r})=\frac{z-z_{n}}{\Delta z}P_{\bm{n}}(\bm{r}),
		\end{gather}
	\end{subequations}
	where $\bm{r}_{\bm{n}}=(x_{n},y_{n},z_{n})$ is the center of each voxel, $P_{\bm{n}}(\bm{r})$ is a volumetric pulse which is equal to $1$ inside voxel $\bm{n}$ and $0$ otherwise. From the above, it follows that the PWL basis functions have support a single element (voxel in our case) of the grid allowing for discontinuities between neighboring voxels. The vector-valued basis function at voxel $\bm{n}$ is:
	\begin{subequations}\label{eq22}
		\begin{gather}
			\bm{f}_{\bm{n}}(\bm{r})=\sum_{q\in\{x,y,z\}}\sum_{l=1}^{4}\bm{f}_{\bm{n}}^{ql}(\bm{r}),\\
			\bm{f}_{\bm{n}}^{ql}(\bm{r})=N_{\bm{n}}^{l}(\bm{r})\hat{\bm{q}},
		\end{gather}
	\end{subequations}
	resulting in $12$ unknowns per voxel. Finally, the electrical properties of the scatterer are modeled by means of PWC approximations:
	\begin{equation}\label{eq23}
		\epsilon_{r}(\bm{r})\approx\sum_{\bm{n}=(1,1,1)}^{(N_{x},N_{y},N_{z})}\epsilon_{r}(\bm{r}_{\bm{n}})P_{\bm{n}}(\bm{r}).
	\end{equation}
	
	\subsection{Galerkin Inner Products and Linear System Formalism}
	
	In what follows, we describe the required steps to numerically solve the current-based formulation of (\ref{eq16}) by means of Galerkin MoM, where the equivalent volumetric currents are expanded in the vector-valued square-integrable basis function (\ref{eq22}) and tested with the same function. At this point, we can form the linear system $\mathbf{A}\mathbf{x}=\mathbf{b}$ with $\mathbf{x},\mathbf{b}\in\mathbb{C}^{N}$ and $\mathbf{A}\in\mathbb{C}^{N\times N}$ with $N=12N_{V}$. More specifically, the matrix and the right-hand side are given by
	\begin{subequations}\label{eq24}
		\begin{gather}
			\mathbf{A}=\mathbf{M}_{\epsilon_{r}}\mathbf{G}-\mathbf{M}_{\chi_e}\mathbf{N},  \\
			\mathbf{b}=c_e\mathbf{M}_{\chi_e}\mathbf{e}_{{\rm inc}},
		\end{gather}
	\end{subequations}
	with
	\begin{subequations}\label{eq25}
		\begin{gather}
			N_{\bm m,\bm n}^{pl,ql'}=\langle\bm{f}_{\bm{m}}^{pl}(\bm{r}),\mathcal{N}\bm{f}_{\bm{n}}^{ql'}(\bm{r})\rangle_{V_{{m}}}, \\
			K_{\bm m,\bm n}^{pl,ql'}=\langle\bm{f}_{\bm{m}}^{pl}(\bm{r}),\mathcal{K}\bm{f}_{\bm{n}}^{ql'}(\bm{r})\rangle_{V_{{m}}},
		\end{gather}
	\end{subequations}
	where $p,q\in\{x,y,z\}$, $l$,$l'\in\{1,2,3,4\}$, $\bm{m}=m_{x}\hat{\mathbf{x}}+m_{y}\hat{\mathbf{y}}+m_{z}\hat{\mathbf{z}}$ denotes the observation voxel and $\bm{n}$ the source voxel. Each voxel interacts with each other so $n_{x},m_{x}=1:N_{x}, n_{y},m_{y}=1:N_{y}$ and $n_{z},m_{z}=1:N_{z}$ resulting in the dense matrices $\mathbf{N},\mathbf{K}\in\mathbb{C}^{N\times N}$, where $\mathcal{K}$ operator is discretized for the computation of the magnetic field, too. Furthermore, $\mathbf{G}\in\mathbb{\mathbb{R}}^{N\times N}$ is the associated Gram matrix, which is diagonal, since non-overlapping basis functions are used:
	\begin{equation}\label{eq26}
		{G}_{\bm m, \bm n}^{pl,ql'}=\langle\bm{f}_{\bm{m}}^{pl}(\bm{r}),\bm{f}_{\bm{n}}^{ql'}(\bm{r})\rangle_{V_{{m}}}.
	\end{equation}
	Also, $\mathbf{M}_{\epsilon_{r}}$, $\mathbf{M}_{\chi_{e}}\in\mathbb{C}^{N\times N}$ are diagonal matrices for isotropic materials with the non-zero values being equal to the material properties at the corresponding voxels. Finally, the "tested" incident electric field that arises from the right-hand side is given by
	\begin{equation}\label{eq27}
		{e}_{{\rm inc}, \bm m}^{pl}=\langle\bm{f}_{\bm{m}}^{pl}(\bm{r}),\bm{e}_{{\rm inc}}(\bm{r})\rangle_{V_{{m}}}.
	\end{equation}
	In the above equations, we use the following definition for the inner products:
	\begin{equation}\label{eq28}
		\langle\bm{f},\bm{g}\rangle_{V}=\intop_{V}{\bar{\bm{f}}}\cdot\bm{g}dV,
	\end{equation}
	where the $\bar{\cdot}$ denotes the complex conjugate operation.
	
	\subsection{Tensor Representation of the Linear System}
	
	The various components of the VIE linear system admit a convenient and intuitive representation in tensor format (multidimensional arrays) when we employ a uniform discretization grid. First, we construct the tensors of the dielectric properties of the scatterer, $\underline{\mathbf{M}}_{\epsilon_{r}}$,$\underline{\mathbf{M}}_{\chi_{e}}\in\mathbb{C}^{N_{x}\times N_{y}\times N_{z}}$, as follows:
	\begin{subequations}\label{eq29}
		\begin{gather}
			\underline{\mathbf{M}}_{\epsilon_{r}}(\bm{m})=\epsilon_{r}(\bm{r}_{\bm{m}}),\\
			\underline{\mathbf{M}}_{\chi_{e}}(\bm{m})=\epsilon_{r}(\bm{r}_{\bm{m}})-1.
		\end{gather}
	\end{subequations}
	Next, we form the tensor of the unknowns, $\underline{\mathbf{x}}^{pl}\in\mathbb{C}^{N_{x}\times N_{y}\times N_{z}}$, where $p\in\{x,y,z\}$ and $l\in\{1,2,3,4\}$. Similarly, we construct the tensor of the incident fields. More specifically, the associated tensor of the  component $p$ and the scalar basis function $l$ reads
	\begin{equation}\label{eq30}
		\underline{\mathbf{e}}_{{\rm inc}}^{pl}(\bm{m})=\intop_{V_{m}}\bm{f}_{\bm{m}}^{pl}(\bm{r})\cdot\bm{e}_{{\rm inc}}(\bm{r})dV=\intop_{V_{m}}N_{\bm{m}}^{l}(\bm{r})e_{{\rm inc}}^{p}(\bm{r})dV,
	\end{equation}
	resulting in $\underline{\mathbf{e}}_{{\rm inc}}^{pl}\in\mathbb{C}^{N_{x}\times N_{y}\times N_{z}}$. Naturally, the tensor of the right-hand side is given by
	\begin{equation}\label{eq31}
		\underline{\mathbf{b}}^{pl}(\bm{m})=c_{e}\underline{\mathbf{M}}_{\chi_{e}}\odot\underline{\mathbf{e}}_{{\rm inc}}^{pl}.
	\end{equation}
	In order to calculate the interactions between the testing and basis functions, we have to calculate the integrals:
	\begin{equation}\label{eq32}
		G_{\bm{m},\bm{n}}^{pl,ql'}=\intop_{V_{m}}\bm{f}_{\bm{m}}^{pl}\cdot\bm{f}_{\bm{n}}^{ql'}dV=\mathbf{\hat{p}}\cdot\mathbf{\hat{q}}\intop_{V_{m}}N_{\bm{m}}^{l}P_{\bm{m}}N_{\bm{n}}^{l'}P_{\bm{n}}dV
	\end{equation}
	for $n_{x}$,$m_{x}=1:N_{x}$, $n_{y}$,$m_{y}=1:N_{y}$ and $n_{z}$,$m_{z}=1:N_{z}$, $p$,$q\in\{x,y,z\}$ and $l$,$l'\in\{1,2,3,4\}$ which can be uniquely expressed by an $\mathbb{R}^{N_{x}\times N_{y}\times N_{z}}$ tensor, since non-overlapping basis functions are used, as
	\begin{equation}\label{eq33}
		\underline{\mathbf{G}}^{pl}(\bm{m})=\intop_{V_{m}}\left( P_{\bm{m}} N_{\bm{m}}^{l}\right)^2dV = \begin{cases}
			\Delta V &l=1\\
			\frac{\Delta V}{12} &l=2,3,4
		\end{cases},
	\end{equation}
	where $\Delta V=\Delta x\Delta y\Delta z$ is the volume of the voxel. 
	
	Finally, the most challenging part in the assembly of the tensors is the calculation of the 6-D integrals that arise from the associated, discretized integro-differential $\mathcal{N}$ and $\mathcal{K}$ operators with kernels that exhibit strongly singular and weakly singular behavior, respectively, when the observation points coincide or are adjacent with the source points. More specifically, these integrals, according to (\ref{eq25}) can be written as
	
	\begin{subequations}\label{eq34}
		\begin{gather}
			N_{\bm{m},\bm{n}}^{pl,ql'}=\intop_{V_{m}}\intop_{V'_{n}}\bm{f}_{\bm{m}}^{pl}(\bm{r})\cdot\nabla\times\nabla\times(g(\bm{r}-\bm{r}')\bm{f}_{\bm{n}}^{ql'}(\bm{r}'))dV'dV,\\		K_{\bm{m},\bm{n}}^{pl,ql'}=\intop_{V_{m}}\intop_{V'_{n}}\bm{f}_{\bm{m}}^{pl}(\bm{r})\cdot\nabla\times(g(\bm{r}-\bm{r}')\bm{f}_{\bm{n}}^{ql'}(\bm{r}'))dV'dV.
		\end{gather}
	\end{subequations}
	
	The far interactions are calculated by means of a standard 6-D quadrature rule, where the components of the dyadic Green function kernels are symmetric and anti-symmetric for $\mathcal{N}$ and $\mathcal{K}$ operator, respectively, resulting in $6$ and $3$ unique dyadic components. Also, the unique interactions of the scalar part of the basis and testing functions with a kernel of the form $K(\bm{r},\bm{r}')=K(\bm{r}-\bm{r}')$ are $10$ instead of $16$ interactions, when discretized on a uniform grid. These internal symmetries significantly reduce the overall memory footprint required for storing the tensors of both operators to $90$ unique entries. However, simple quadrature rules cannot be applied when calculating the nearby interactions, since for coinciding or adjacent voxels the 6-D integrals are singular. The main idea for tackling these cases is to reduce the dimensionality of the volume-volume integrals to surface-surface step-by-step, as it has been shown in \cite{georgakis2017reduction,knockaert2011analytic}. As a result, the initial volume-volume integral boils down to a sum of $4$  surface-surface integrals, over the faces of the interacting voxels, which have smoother kernels and can be calculated by modern algorithms, originally developed for Galerkin SIE methods over quadrilateral patches \cite{DIRECTFN,tambova2018generalization}.  
	
	In this work, we choose to use a uniform, voxelized grid which enables us to exploit the translation invariance property of the convolutional discrete kernels:
	\begin{equation}\label{eq42}
		N_{\bm{m},\bm{n}}^{pl,ql'} = N_{\bm{m}-\bm{n}}^{pl,ql'},
	\end{equation}
	which means that we can fix the basis function at a specific voxel, say the first with $\bm{n}=(1,1,1)$, and sweep the testing function over the voxels of the computational domain in order to calculate only the unique volume-volume integrals. Therefore, the above integrals, when expressed in matrix $\in\mathbb{C}^{N_{V}\times N_{V}}$ format will be three-level block-Toeplitz Toeplitz-block matrices, since we are dealing with the 3-D case. As a result, it suffices to calculate the Toeplitz defining tensors $\in\mathbb{C}^{N_{x}\times N_{y}\times N_{z}}$, as
	\begin{equation}\label{eq43}
		\underline{\mathbf{N}}^{pl,ql'}(\bm{m})=N_{\bm{m}-\bm{1}}^{pl,ql'},
	\end{equation}
	where $m_{x}=1:N_{x}$, $m_{y}=1:N_{y}$, $m_{z}=1:N_{z}$, $p$,$q\in\{x,y,z\}$ and $l$,$l'\in\{1,2,3,4\}$. The same logic applies to the tensors for $\mathcal{K}$ operator. Finally, even greater compression can be achieved by applying low multilinear rank tensor decompositions on the Toeplitz defining tensors \cite{polimeridis2014compression,giannakopoulos20183D,giannakopoulos2019memory}.
	
	\section{FFT-BASED VIE SOLVER}
	
	\subsection{Acceleration of the Matrix-Vector Product}
	
	In order to numerically solve a discretized volume integral equation, the inversion of a very large and dense matrix is required. Hence, the scientific community has naturally been using iterative solvers. However, even with an iterative solver, for which the most time consuming part is the matrix-vector product with $\mathcal{O}(N^{2})$ complexity, solving a realistic problem becomes extremely difficult and in many cases prohibitive. However, if uniform, voxelized grids are used, the MoM matrix becomes structured (block Toeplitz with Toeplitz blocks) and the associated matrix-vector product is amenable to acceleration with the help of FFT resulting in $\mathcal{O}(N\text{log}(N))$ complexity for each matrix-vector product. 
	
	The first step to that direction is the embedding of the Toeplitz defining tensors for $\mathcal{N}$ and $\mathcal{K}$ operators into their respective circulant defining tensors $\underline{\mathbf{N}}_{{\rm circ}}^{pl,ql'}\in\mathbb{C}^{2N_{x}\times2N_{y}\times2N_{z}}$, with special care to take into account the associated projection signs for the dyadic Green functions, the basis, and testing functions, as described in a  detailed algorithmic analysis, that can be found at the Appendix of \cite{polimeridis2014stable}, which we omit here for brevity. Finally, the circulant matrix-vector product can be computed with the FFT, reducing the computational cost to approximately $\mathcal{O}(N\text{log}(N))$.
	
	\subsection{On Preconditioning}
	Furthermore, an advantage of using the current-based formulation is that the multiplicative operator with the dielectric susceptibilities is on the left of $ \mathcal{N} $ operator, as it can be seen in (\ref{eq16}). Therefore, by multiplying through by the inverse of $\mathcal{M}_{\epsilon_r}$, this equation can be rewritten as 
	\begin{equation}\label{eq17}
		(\mathcal{I}-\mathcal{M}_{\tau_{e}}\mathcal{N})\bm{j}=c_{e}\mathcal{M}_{\tau_{e}}\bm{e}_{{\rm inc}},
	\end{equation}
	with $\tau_{e}\triangleq\chi_{e}/\epsilon_{r}$, where the identity operator is now left alone. In this way, the integral equation is regularized and this regularization can be thought as a natural preconditioning for solving the linear system of the discretized current-based VIE, as it has been suggested in \cite{zouros2012transverse,budko2006spectrum}. Motivated from these observations, we use a preconditioner of the form $\mathbf{P}=\mathbf{M}_{\epsilon_{r}}\mathbf{G}$. From a numerical perspective, when using this preconditioner, the iterative solver converges much faster especially in the case of highly inhomogeneous scatterers. Numerical experiments, presented below, indicate the superior convergence properties of the iterative solver when using this preconditioned solver for the case of a RHBM.
	
	\subsection{Computation of Fields and Power}
	The linear system solution provides us with the equivalent electric and magnetic currents at the center of the voxels and the three linear coefficients represent the slopes of the linear terms. However, these quantities are mathematical inventions and do not represent any physical phenomenon or quantity. In this subsection, we describe the required steps for the computation of the EM fields and the associated absorbed power. 
	
	A fast and efficient way to calculate these quantities is by employing the fast matrix-vector products and calculating the "tested" total fields according to the discrete analogue of (\ref{tot_fields}). Furthermore, it should be noted that, in order to calculate the fields in a discrete sense, we need to divide the "tested" fields by the Gram matrix. Finally, an advantage of this approach is that no extra memory is required for computing the fields, and the complexity is once again governed by the FFT-based matrix-vector product.
	
	Finally, the absorbed power can be accurately computed based on \cite{polimeridis2015computation}, where stable and compact vector-matrix-vector formulas are suggested. The only difference here is that we have to replace the Gram matrix, as in (\ref{eq33}), to include the contributions of the linear terms of the solution vector and compute the absorbed power, which corresponds also to the SAR, as
	\begin{equation}\label{eq44}
		P_{\rm abs} = \frac{1}{2}\textrm {Re} \left\{ \mathbf{x}^* ( c_e \mathbf{M}_{\chi_e} )^{-1} \mathbf{G} \mathbf{x} \right\},
	\end{equation}
	with the $^*$ superscript denoting the conjugate transpose operation.
	
	\section{NUMERICAL RESULTS}
	
	In this Section, we report results for representative numerical experiments that illustrate the convergence properties and numerical accuracy of the proposed solver. Specifically, we validate the proposed solver by comparison with the Mie series solution for dielectric spheres, study its convergence rate for dielectric cubes, and demonstrate its accuracy and the convergence properties of the iterative solver when modeling the EM scattering from RHBMs. For all PWC, PWL, and DVIE computations the FFT is employed and the iterative solver of choice is the Generalized Minimum Residual (GMRES) with inner iterations 50, outer iterations 200, and tolerance $10^{-5}$.
		
	The computations were performed with Matlab (version 9.2) at a CentOS 6.9 operating system equipped with an Intel(R) Xeon(R) central processing unit (CPU) E5-2699 v3 @ 2.30 GHz with 36 cores and 2 threads per core and an NVIDIA Tesla K40M GPU. We report the computational time that GMRES requires to converge to the specified tolerance, as well as, the memory footprint of the tensors of the unknowns for PWC, PWL, and DVIE. These tensors have dimensions of $2N_{x}\times2N_{y}\times2N_{z} = 8N_V$ and the number of unknowns per voxel are 12 for PWL, 3 for PWC, and 3 for DVIE, resulting in an unknowns' total memory footprint of $96N_V$, $24N_V$, and $24N_V$ for PWL, PWC, and DVIE respectively. For completeness, we note that each circulant defining tensor has dimensions of $2N_{x}\times2N_{y}\times2N_{z}$ and $\mathcal{N}$ operator consists of 60 unique tensors for PWL, 6 for PWC, and the DVIE Green function requires only 1 tensor. However, such tensors are compressible \cite{giannakopoulos2019memory}, therefore, we do not include the memory they require when we compute the memory footprint of the solvers.
	
	\subsection{Convergence Study for Dielectric Spheres}
	
	To validate the proposed framework quantitatively, we conduct the numerical computation for the scattering of a plane wave from homogeneous dielectric spheres and calculate the absorbed power, as in (\ref{eq44}). The results are compared with the Mie series solution \cite{mie1908beitrage}. In detail, we model frequency-dependent electrical properties, matched to the average of the gray and white matter according to \cite{gabriel1996compilation,guerin2017ultimate}, as they can be seen in Table \ref{dielectric_properties}, for various $B_0$ static magnetic field strengths of an MR scanner, where we use the gyromagnetic ratio as $\frac{\gamma}{2\pi}=42.58$ $\rm MHz/T$ to set the operating frequency $f=\frac{\gamma}{2\pi}B_0$. Furthermore, the sphere has radius $r=7.5$ $\rm cm$ and the excitation is an $x$-polarized and $z$-directed plane wave $\bm{e}_{\rm inc} = \bm{\hat{x}}e^{-ik_0z}$.
	
	\begin{table}[ht!]
		\captionsetup{font=scriptsize}
		\caption{ELECTRICAL PROPERTIES} \label{dielectric_properties} \centering
		\begin{tabular}{c|cccc}
			
			\hline\hline\\[-0.4em]
			$B_0(\rm T)$    		& 0.5   & 1.5   & 3.0      & 7.0       \\[0.3em] \hline \\[-0.3em]
			$f(\rm MHz)$	        & 21.29 & 63.87 & 127.74 & 298.06  \\[0.3em]
			$\epsilon_r^{'}$		& 227   & 141   & 68     & 64	   \\[0.3em]
			$\sigma(\rm S/m)$		& 0.25  & 0.35  & 0.50    & 0.51	   \\[0.3em]
			
			\hline\hline
		\end{tabular}
	\end{table}
	
	\begin{table}[ht!]
		\captionsetup{font=scriptsize}
		\caption{DoFs AND MEMORY FOR EXPERIMENT A} \label{dofs_A} \centering
		\begin{tabular}{c|ccccc}
			\hline\hline\\[-0.4em]
			$h(\rm mm)$    		        & 5.00      & 2.50     & 1.50     & 1.00       & 0.75     \\[0.3em] \hline \\[-0.3em]
			$N_x=N_y=N_z$				& 31     & 61      & 101     & 151     & 201      \\[0.3em]
			DoFs PWC (million)			& 0.09   & 0.68    & 3.09    & 10.33   & 24.36    \\[0.3em]
			DoFs PWL (million)	        & 0.36   & 2.72    & 12.36   & 41.32   & 97.45    \\[0.3em]\hline&&&&\\[-0.4em]
			PWC Unknowns Mem. (GB)		& 0.011  & 0.081   & 0.368   & 1.231   & 2.904    \\[0.3em]
			PWL Unknowns Mem. (GB)	    & 0.043  & 0.325   & 1.474   & 4.925   & 11.617   \\[0.3em]
			\hline\hline
		\end{tabular}
	\end{table}
	
	\begin{table}[ht!]
		\captionsetup{font=scriptsize}
		\caption{COMPUTATIONAL TIMES FOR EXPERIMENT A} \label{times_A} \centering
		\begin{tabular}{c|c|cccc}\hline\hline&&&&\\[-0.4em]
			Solver  &  \backslashbox{${h}$}{${B_0}$}    & 0.5 T    & 1.5 T   & 3.0 T       & 7.0 T \\[0.3em] \hline &&&&\\[-0.3em]
			&    5.00 mm                      & 00:00:05  & 00:00:02   & 00:00:02  & 00:00:02      \\[0.3em]
			&    2.50 mm 	    		  	  & 00:00:14  & 00:00:11   & 00:00:10  & 00:00:13     \\[0.3em]
   PWC      &	 1.50 mm 					  & 00:01:10  & 00:00:58   & 00:00:53  & 00:01:07     \\[0.3em]
			&	 1.00 mm       				  & 00:04:15  & 00:03:37   & 00:03:18  & 00:04:10    \\[0.3em]
			&	 0.75 mm	    			  & 00:16:14  & 00:13:43   & 00:12:28  & 00:16:06    \\[0.3em] \hline&&&&\\[-0.4em]
			&    5.00 mm                      & 00:00:08  & 00:00:06   & 00:00:05  & 00:00:08      \\[0.3em]
			&    2.50 mm 		    		  & 00:00:34  & 00:00:29   & 00:00:23  & 00:00:37     \\[0.3em]
   PWL		&	 1.50 mm  				  & 00:33:10  & 00:29:50   & 00:26:16  & 00:36:28     \\[0.3em]
			&	 1.00 mm	       			  & 01:43:38  & 01:32:55   & 01:22:54  & 01:56:10   \\[0.3em]
			&	 0.75 mm 	    		  & 03:59:05  & 03:29:05   & 03:09:38  & 04:36:17  \\[0.3em]
			\hline\hline
		\end{tabular}
	\end{table}

\begin{figure}[ht!]
		\captionsetup{font=footnotesize}
		\centering
		\includegraphics[scale=0.43]{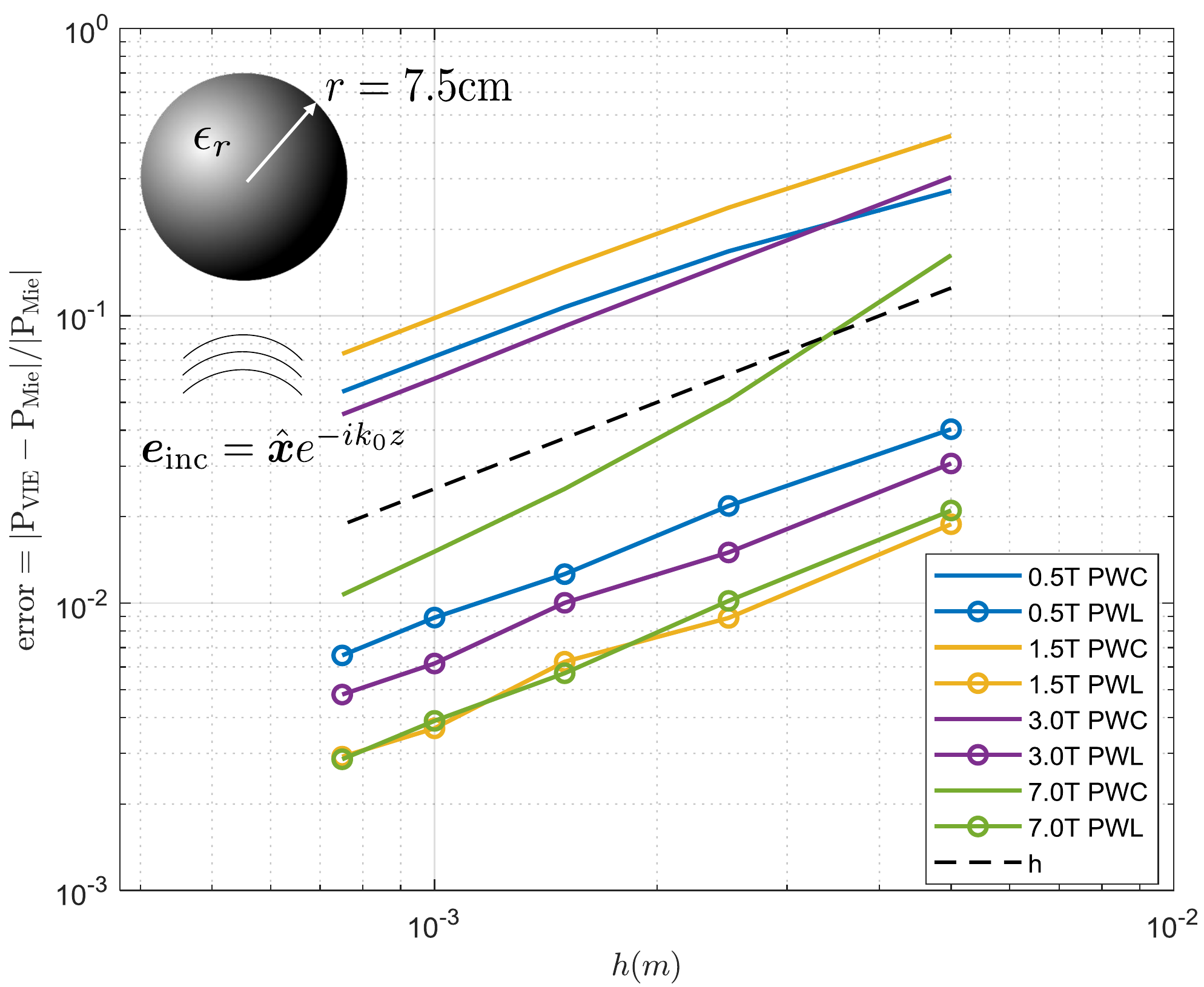}
		\caption{Relative error of the absorbed power, with respect to the resolution, for the PWC and PWL current-based VIE, compared with the Mie series solution, for homogeneous spheres with $r = 7.5$ $\rm cm$ and frequency-dependent electrical properties, as in Table \ref{dielectric_properties}, for each different $B_0$, when irradiated by a plane wave.}
		\label{sphere_refinements}
	\end{figure}

	Moreover, we calculate the absorbed power, both with the PWC current-based VIE \cite{polimeridis2014stable} and the proposed solver with PWL basis functions, and compute the relative error for the power obtained from Mie series, as $\rm error = |P_{\rm VIE} - P_{\rm Mie}| / |P_{\rm Mie}|$. Finally, we vary the resolution of the discretization, using $h = 5$, $2.5$, $1.5$, $1$, and $0.75$ $\rm mm$, and present the results for the relative error of the PWC and PWL solvers, with respect to the resolution in Fig. \ref{sphere_refinements}. The Degrees of Freedom (DoFs) and memory footprint of each solver are shown in Table \ref{dofs_A} and Table \ref{times_A} shows the computational time (in [h]:mm:ss) that GMRES requires to converge. It should be noted that the matrix-vector product is performed in the GPU up to 1 mm resolution for PWC and up to 2.5 mm for PWL, thus, the significant increase in the computational time.
	
	Clearly, the discretization with PWL basis functions gives much better accuracy with respect to the resolution, compared with the PWC discretization. Specifically, it is observed that there is a constant improvement factor of approximately 12 times in the relative error between the PWC and PWL solvers. However, the convergence rate appears to be $\mathcal{O}(h)$, instead of $\mathcal{O}(h^2)$, due to numerical inaccuracies originating from the staircase approximation of a sphere, when discretized with voxels.

	\subsection{Convergence Study for Dielectric Cubes}
	
	To overcome the numerical inaccuracies associated with the staircase approximation of a spherical object in a voxel-based discretization, in this subsection, we compute the scattering of a plane wave from homogeneous dielectric cubes instead of spheres and calculate the absorbed power. The results are compared with the SIE open-source solver \textit{scuff-em} \cite{SCUFF}, where 24,576 edge elements have been used to ensure high accuracy in the reference solution. The cube has an edge of $L = 15$ $\rm cm$ and the rest parameters, including the electrical properties, the operating frequencies, and the excitation, are identical to those of Section V.A (see Table \ref{dielectric_properties}).
	
	\begin{table}[ht!]
		\captionsetup{font=scriptsize}
		\caption{DoFs AND MEMORY FOR EXPERIMENT B} \label{dofs_B} \centering
		\begin{tabular}{c|ccccc}
			\hline\hline\\[-0.4em]
			$h(\rm mm)$  		& 5.00       & 2.50      & 1.50     & 1.00       & 0.75     \\[0.3em] \hline \\[-0.3em]
			$N_x=N_y=N_z$		& 30      & 60       & 100     & 150     & 200      \\[0.3em]
			DoFs PWC (million)	& 0.081   & 0.648    & 3       & 10.125  & 24       \\[0.3em]
			DoFs PWL (million)	& 0.324   & 2.592    & 12      & 40.500  & 96       \\[0.3em] 
			DoFs DVIE (million)	& 0.081   & 0.648    & 3       & 10.125  & 24       \\[0.3em]\hline&&&&\\[-0.4em]
			PWC Unknowns Mem. (GB)		& 0.010    & 0.077    & 0.358   & 1.207   & 2.861    \\[0.3em]
			PWL Unknowns Mem. (GB)	    & 0.039    & 0.309    & 1.431   & 4.828   & 11.444   \\[0.3em]
			DVIE Unknowns Mem. (GB)	    & 0.010    & 0.077    & 0.358   & 1.207   & 2.861    \\[0.3em]
			\hline\hline
		\end{tabular}
	\end{table}
	
	\begin{table}[ht!]
		\captionsetup{font=scriptsize}
		\caption{COMPUTATIONAL TIMES FOR EXPERIMENT B} \label{times_B} \centering
		\begin{tabular}{c|c|cccc}\hline\hline&&&&\\[-0.4em]
			Solver  &  \backslashbox{${h}$}{${B_0}$}    & 0.5 T    & 1.5 T   & 3.0 T       & 7.0 T \\[0.3em] \hline &&&&\\[-0.3em]
			&    5.00 mm                  & 00:00:03  & 00:00:02   & 00:00:02  & 00:00:04      \\[0.3em]
			&    2.50 mm 	    		  & 00:00:18  & 00:00:16   & 00:00:16  & 00:00:32     \\[0.3em]
	PWC		&	 1.50 mm 				  & 00:01:20  & 00:01:12   & 00:01:11  & 00:02:27     \\[0.3em]
			&	 1.00 mm       			  & 00:04:36  & 00:04:09   & 00:04:07  & 00:08:25    \\[0.3em]
			&	 0.75 mm	    		  & 00:18:03  & 00:16:02   & 00:15:50  & 00:32:09    \\[0.3em] \hline&&&&\\[-0.4em]
			
			&    5.00 mm                  & 00:00:09  & 00:00:08   & 00:00:08  & 00:00:15      \\[0.3em]
			&    2.50 mm 	    	 	  & 00:01:18  & 00:01:10   & 00:01:08  & 00:02:22     \\[0.3em]
	PWL 	&	 1.50 mm  			  & 00:30:23  & 00:26:32   & 00:26:11  & 00:52:37     \\[0.3em]
			&	 1.00 mm	       		  & 01:38:47  & 01:26:40   & 01:25:02  & 02:52:51   \\[0.3em]
			&	 0.75 mm 	    	  & 03:47:41  & 03:16:53   & 03:14:24  & 06:43:00  \\[0.3em]\hline&&&&\\[-0.4em]
			
			&    5.00 mm                  & 00:00:16  & 00:00:28   & 00:00:27  & 00:01:09      \\[0.3em]
			&    2.50 mm 		    	  & 00:01:07  & 00:01:49   & 00:01:49  & 00:05:11     \\[0.3em]
	DVIE	&	 1.50 mm  			  & 00:12:22  & 00:15:32   & 00:06:41  & 00:46:26     \\[0.3em]
			&	 1.00 mm	       		  & 00:41:11  & 00:51:09   & 01:12:20  & 02:53:50   \\[0.3em]
			&	 0.75 mm 	    	  & 01:22:00  & 01:56:33   & 02:15:27  & 06:01:55  \\[0.3em]
			\hline\hline
		\end{tabular}
	\end{table}
	
	\begin{figure}[ht!]
		\captionsetup{font=footnotesize}
		\centering
		\includegraphics[scale=0.43]{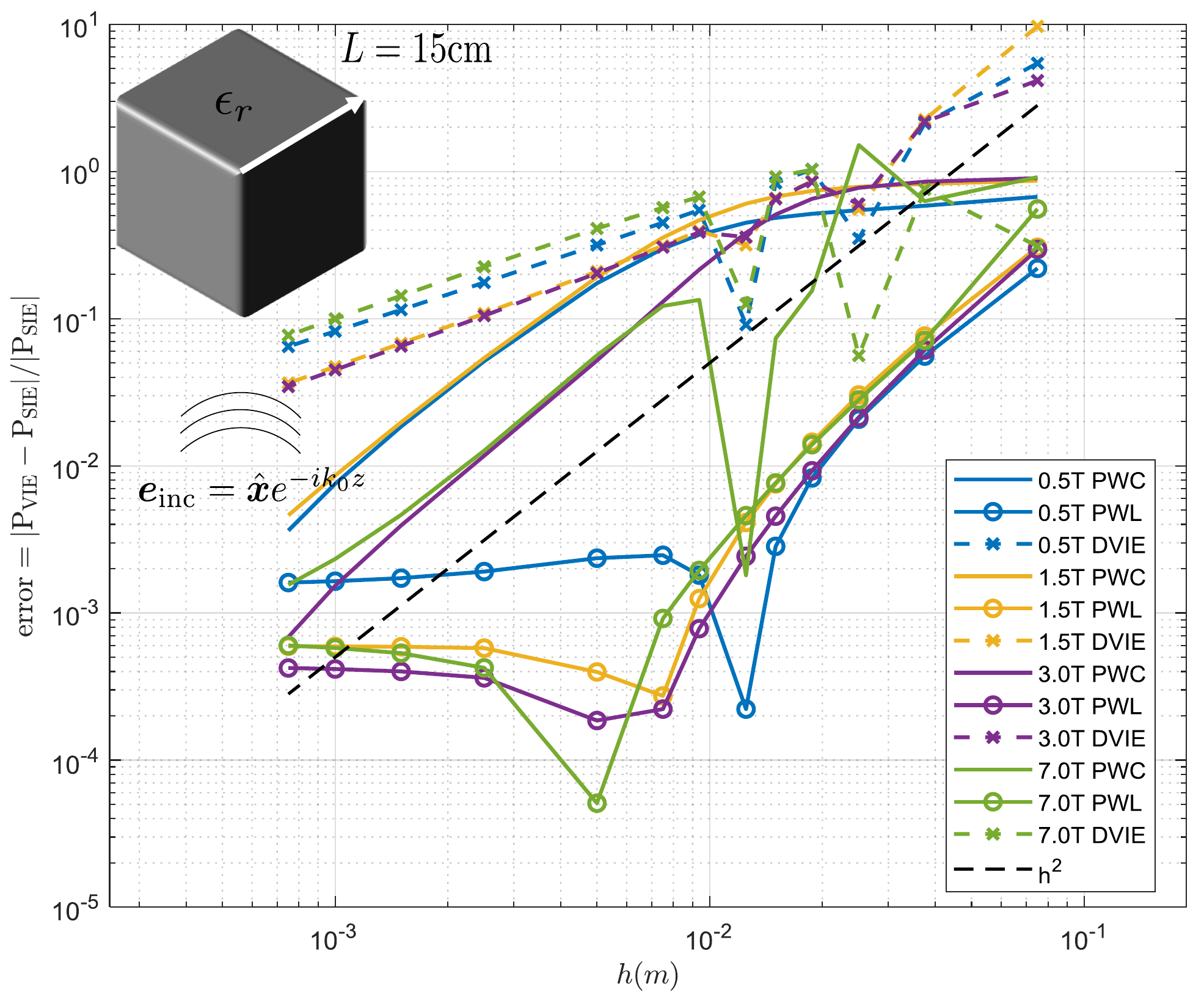}
		\caption{Relative error of the absorbed power, with respect to the resolution, for the PWC and PWL current-based VIE and the flux-based DVIE solver compared with the SIE solution, for homogeneous cubes with $L = 15$ $\rm cm$ and frequency-dependent electrical properties, as in Table \ref{dielectric_properties}, for each different $B_0$, when irradiated by a plane wave.}
		\label{cube_refinements}
	\end{figure}
		
	We calculate the absorbed power, both with the PWC current-based VIE \cite{polimeridis2014stable}, the proposed PWL solver, and with the flux-based \cite{zwamborn1992three,zwanborn1991weak} (DVIE) solver. We then compute the relative error for the power obtained from \textit{scuff-em}, as $\rm error = |P_{\rm VIE} - P_{\rm SIE}| / |P_{\rm SIE}|$. This time, we vary the resolution of the voxels in a broader range, using $h = 75$, $37.5$, $25$, $18.75$, $15$, $12.5$, $9.375$, $7.5$, $5$, $2.5$, $1.5$, $1$, and $0.75$ $\rm mm$, and present the results for the relative error of the PWC, PWL, and DVIE solvers, with respect to the resolution in Fig. \ref{cube_refinements}. The DoFs and memory footprint of each solver are shown in Table \ref{dofs_B} and the computational time that GMRES requires to converge in Table \ref{times_B}. We show the computational resources for resolutions finer than 5 mm, adding consistency with Section V.A and brevity, since for coarser resolutions the maximum memory and time required is only 70MB and 8 seconds, respectively. It should also be noted that the matrix-vector product is performed in the GPU up to 1 mm resolution for PWC and DVIE and up to 2.5 mm for PWL, as in the case of the sphere.
	
	Similarly to the sphere case, PWL provides much better accuracy than both PWC and DVIE solvers, especially at coarser resolutions and higher field strengths. It is interesting to note that at fine resolutions PWL relative error reaches a plateau, which can be explained by the fact that the average edge length of the SIE-based reference solution is $4.5$ mm. Therefore, such parameter sets the best spatial resolution that can be achieved by the reference solution and reducing the resolution of the VIE PWL solver below this, would not improve the accuracy. However, for coarser resolutions PWL solver converges, as expected, with an $\mathcal{O}(h^2)$ rate while for the same coarse resolutions PWC has not started to converge. Finally, we note that DVIE solver converges with an $\mathcal{O}(h)$ rate.

	\subsection{Convergence Study for the Current-Based and Flux-Based Formulation}
	
	In this subsection, in order to demonstrate the more stable convergence properties of the proposed high-order current-based solver, in comparison with the flux-based formulation \cite{zwamborn1992three,zwanborn1991weak}, often used in MRI studies \cite{jin1996computation,huang2018fast}, we conduct the computations for the EM scattering of a plane wave for the "Billie" RHBM from Virtual Family Population \cite{christ2009virtual}. In detail, we operate at $B_0 = 7$ $\rm T$, hence the corresponding frequency is $f=298.06$ $\rm MHz$ and we irradiate the RHBM with an $x$-polarized and negative $y$-directed plane wave $\bm{e}_{\rm inc} = \bm{\hat{x}}e^{ik_0y}$, so that the plane wavefront is parallel to the coronal plane of the RHBM.
	
\begin{figure}[ht!]
		\captionsetup{font=footnotesize}
		\centering
		\includegraphics[scale=0.9]{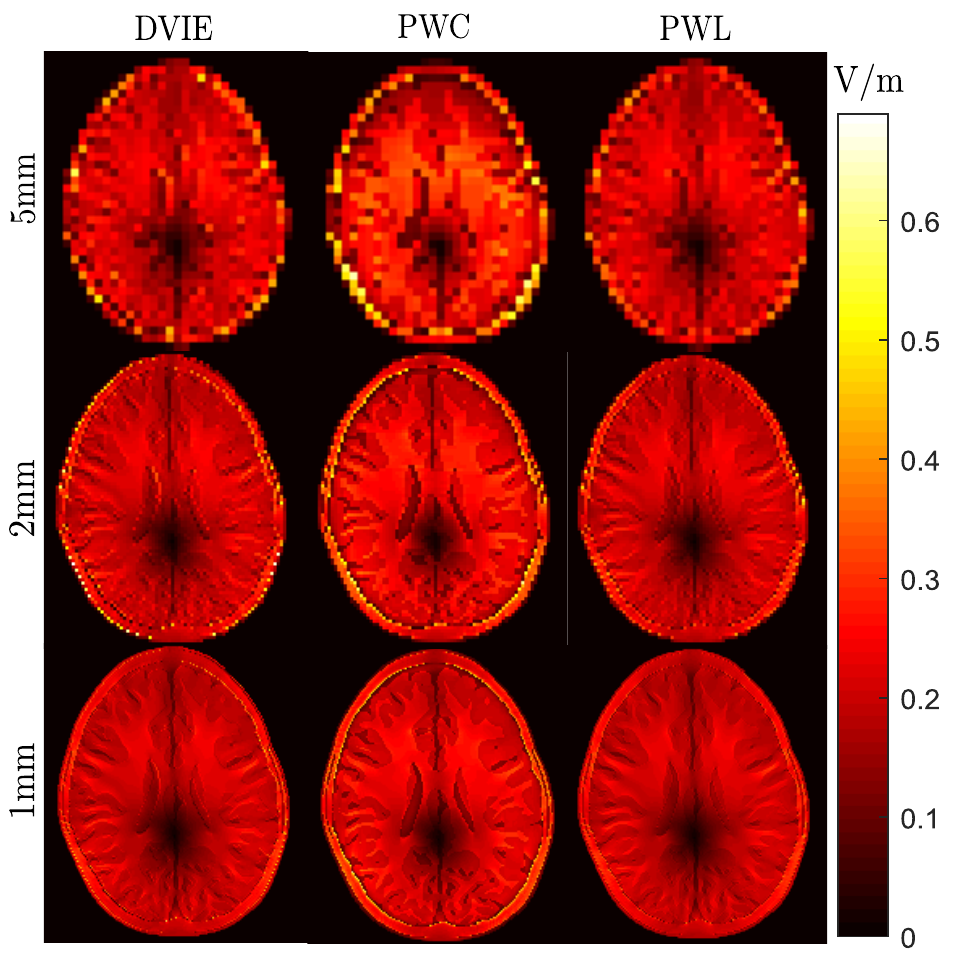}
		\caption{Axial views of the root mean square (RMS) value of the electric field in the "Billie" RHBM, when excited by a plane wave at $7$ $\rm T$ MRI. From left to right, the electric field values are obtained with the flux-based solver (DVIE), the PWC current-based solver, and the PWL current-based solver for resolutions  $h = 5$, $2$, and $1$ $\rm mm$, as viewed from top to bottom. Fields outside the body are masked to improve the visibility.}
		\label{jvie_vs_dvie_etot}
		
	\end{figure}

	Subsequently, we compute the EM fields and the absorbed power with the current-based formulation with PWC and PWL basis functions and the flux-based (DVIE) formulation with rooftop basis functions, as we refine the grid, using resolutions $h = 5$, $2$, and $1$ $\rm mm$. Table \ref{dofs_times_C} shows the DoFs and the memory footprint of each solver and the computational time that the iterative solver requires to converge for each method. The matrix-vector product is performed in the GPU up to 2 mm for PWC and DVIE and up to 5 mm for PWL. In Fig. \ref{jvie_vs_dvie_etot}, we present the root mean square (RMS) electric field at an axial cut for all $3$ resolutions and $3$ solvers, where the fields outside the body are masked for enhanced visualization. 
	
	At Table \ref{dofs_times_C}, we also show the absorbed power inside the head for different solvers and resolutions, where we notice a decreasing trend with $h$-refinement for all solvers. Additionally, PWL absorbed powers are lower than PWC ones for the same resolution and 5 mm PWL has a lower power than 2 mm PWC as well as 2 mm PWL than 1 mm PWC. Therefore, under the assumption that a lower absorbed power implies a lower error, it could be argued that coarser PWL electric fields are more accurate than finer PWC electric fields. Additionally, to perform a completely quantitative analysis, we consider the finest 1 mm PWL electric fields as the reference solution and evaluate the following errors:
	\begin{subequations}
		\begin{gather}
			 \textrm{error}_{\bm e_1} = |  P_{\rm abs} - \hat{P}_{\rm abs} |    / | \hat{P}_{\rm abs} |, 		\\
			 \textrm{error}_{\bm e_2} = || \mathbf p_{\rm abs} - \hat{\mathbf p}_{\rm abs} ||_1 / ||\hat{\mathbf p}_{\rm abs} ||_1, 	\\  
			 \textrm{error}_{\bm e_3} = || \mathbf e   - \hat{\mathbf e}	||_2 / ||\hat{\mathbf e}   ||_2, 
		\end{gather}
	\end{subequations}
	where the quantities under $\hat{\cdot}$ are the reference quantities and $\mathbf p_{\rm abs} = 1/2\sigma|\mathbf e|^2$ is the absorbed power density. At Table \ref{dofs_times_C}, we show these errors and notice that PWL 5 mm is more accurate than PWC 2 mm for all these errors despite PWC having $\sim 3.75$ times the unknowns of PWL. For the comparison between PWL 2 mm and PWC 1 mm, we see that PWL $\textrm{error}_{\bm e_1}$ is smaller than PWC $\textrm{error}_{\bm e_1}$, $\textrm{error}_{\bm e_2}$  is almost equal for PWC and PWL, and PWL $\textrm{error}_{\bm e_3}$ is larger than PWC $\textrm{error}_{\bm e_3}$, despite PWC having twice the DoFs of PWL. Finally, for the case of the PWC solver there exist numerical artifacts at regions with high contrast, i.e., at the air-skull and skull-brain interfaces, which is not the case for the proposed higher-order solver. 
	
	\begingroup
	\setlength{\tabcolsep}{3pt}
	\begin{table}[ht!]
		\captionsetup{font=scriptsize}
		\caption{COMPUTATIONAL CONSIDERATIONS FOR EXPERIMENT C}  \label{dofs_times_C} \centering
		\begin{tabular}{c|ccc}
			\hline\hline\\[-0.4em]
			$h(\rm mm)$ 		        & 5      & 2    & 1            \\[0.3em] \hline \\[-0.3em]
			$N_x\times N_y\times N_z$	& 34$\times$38$\times$45   & 84$\times$94$\times$111    & 168$\times$188$\times$222 \\[0.3em]
			DoFs DVIE (million)         & 0.17   	 & 2.63        & 21.03     \\[0.3em]
			DoFs PWC (million)			& 0.17   	 & 2.63        & 21.03     \\[0.3em]
			DoFs PWL (million)	        & 0.70    	 & 10.52       & 84.14     \\[0.3em]\hline&&&\\[-0.4em]
			DVIE Unknowns Mem. (GB)	    	    & 0.021  	 & 0.313       &  2.508     \\[0.3em]
			PWC Unknowns Mem. (GB)				& 0.021  	 & 0.313       &  2.508      \\[0.3em]
			PWL Unknowns Mem. (GB)	    		& 0.083  	 & 1.254       & 10.030    \\[0.3em]\hline&&&\\[-0.4em]
			DVIE Time                   & 00:00:16   & 00:02:55    & 01:53:33     \\[0.3em]
			PWC Time w/ prec			& 00:00:31   & 00:02:02    & 00:52:41     \\[0.3em]
			PWC Time w/o prec	        & 00:00:57   & 00:04:12    & 01:42:19     \\[0.3em]
			PWL Time w/ prec	    	& 00:02:09   & 01:30:12    & 22:37:37    \\[0.3em]
			PWL Time w/o prec			& 00:11:16   & 02:42:00    & 102:11:12   \\[0.3em]\hline&&&\\[-0.4em]
			DVIE $P_{\rm abs}(\rm \mu W)$	      	 & 29.03  	   & 27.63    & 27.27     \\[0.3em]
			PWC  $P_{\rm abs}(\rm \mu W)$			 & 30.77       & 28.50    & 27.68      \\[0.3em]
			PWL  $P_{\rm abs}(\rm \mu W)$	    	 & 27.82       & 27.00    & 26.65    \\[0.3em]\hline&&&\\[-0.4em]
			PWC $\textrm{error}_{\bm e_1}$ (\%) 				&15.43	 & 6.94	   		&3.86		\\[0.3em]
			PWL $\textrm{error}_{\bm e_1}$ (\%)					& 4.39	 & 1.31	   		&N/A			\\[0.3em]
			PWC $\textrm{error}_{\bm e_2}$ (\%)					&63.18	 &41.80	   		&26.19		\\[0.3em]
			PWL $\textrm{error}_{\bm e_2}$ (\%)					&38.86	 &27.49	        &N/A			\\[0.3em]
			PWC $\textrm{error}_{\bm e_3}$ (\%)					&69.68	 &51.43	        &25.05		\\[0.3em]
			PWL $\textrm{error}_{\bm e_3}$ (\%)					&43.92	 &38.12	        &N/A			\\[0.3em]
			\hline\hline
		\end{tabular}
	\end{table}
	\endgroup
	
	\begin{figure}[ht!]
		\captionsetup{font=footnotesize}
		\centering
		\includegraphics[scale=0.3]{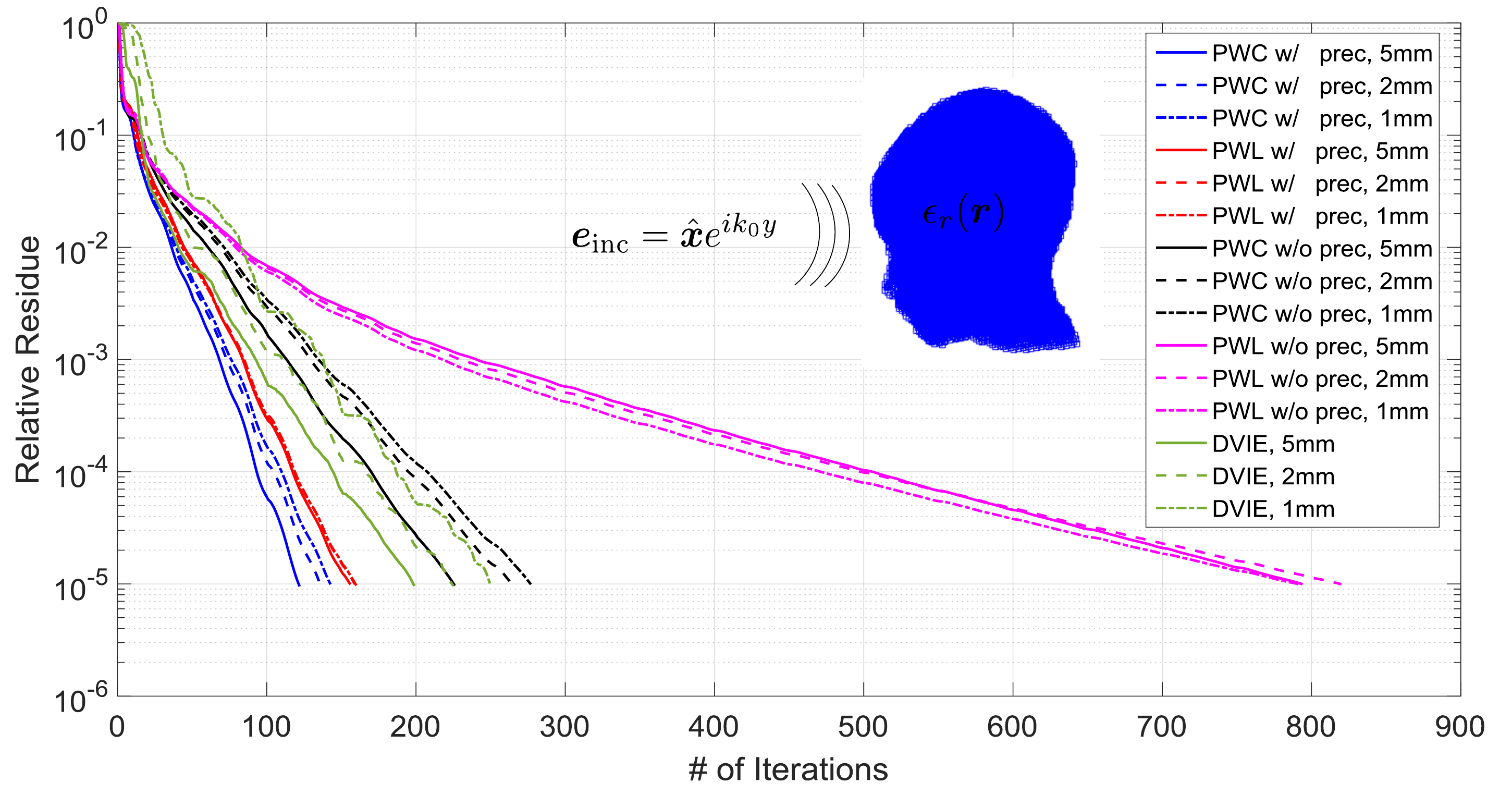}
		\caption{Convergence of GMRES iterative solver for the PWC current-based solver, the PWL current-based solver with and without preconditioner, and the flux-based solver (DVIE) for different resolutions $h = 5$, $2$, and $1$ $\rm mm$, for the calculation of the EM scattering of a plane wave from the "Billie" RHBM at $7$ $\rm T$ MRI.}
		\label{convergence-jvie-dvie}
	\end{figure}
	
	Moreover, we execute the same numerical experiment, utilizing the preconditioner mentioned in the previous section, and present the convergence of the iterative solver for all formulations and refinements. As it can be clearly seen in Fig. \ref{convergence-jvie-dvie}, the current-based formulation is well-conditioned, allowing for $p$-refinement, while the iteration count remains practically the same, when the suggested preconditioner is used. It also allows for $h$-refinement, since, as we move to finer resolutions, the number of iterations does not change. On the contrary, the iteration count of the flux-based formulation in this particular example demonstrates a stronger dependence on the number of unknowns and in the limiting case, when the relative permittivity goes to infinity, DVIE becomes a first-kind integral equation \cite{markkanen2012analysis,polimeridis2014stable} with its associated convergence issues with grid refinement. Also, as shown in \cite{yla2014surface} for a dielectric sphere, the DVIE iteration count is expected to be more dependent on the number of elements than the one from the current-based VIE for which the spectral properties of the operators are preserved. The iteration count increase could also be associated to the fact that the flux-based VIE, in \cite{zwamborn1992three}, averages the normalized contrast function with the material properties.
	
	The convergence of the unpreconditioned PWL solver is significantly slower than the corresponding of the PWC solver, which should be associated with the $p$-refinement itself and with the fact that the PWC basis functions are orthogonal while the PWL are not. Such property leads to an increased condition number of the PWL Gram matrix and, subsequently, of the unpreconditioned PWL system matrix. However, the proposed preconditioner, without practically adding any computational overhead, drops significantly the iteration count, particularly for the case of the PWL solver, demonstrating its effectiveness for highly inhomogeneous scatterers. Finally, it should be noted that such preconditioner can only be applied in the current-based VIE and is not available for the flux-based formulation.
	
	\subsection{Convergence Study for the Case of a High-Permittivity Pad Attached to the RHBM and Comparison with FDTD}
	
	In MRI applications, the transmit magnetic field defined as $b_1^+ = \mu_0 (h_x + ih_y)$ is required to be homogeneous to uniformly tip the magnetization vector. In high-field MRI, where the effective wavelength is comparable to the dimensions of the subject, $b_1^+$ becomes inhomogeneous and the use of high-permittivity dielectric materials has been proposed \cite{webb2011dielectric,de2012increasing,brink2014high,brink2015effect,koolstra2017improved,haemer2018approaching,vaidya2018manipulating,vaidya2018improved} as an effective way to address such $b_1^+$ inhomogeneities. By introducing specifically designed dielectric pads between the RF coil and the subject, $b_1^+$ becomes more homogeneous, SAR is reduced and SNR increases, a technique called dielectric shimming. However, from an EM scattering perspective, simulating such setups is very challenging even for the state-of-the-art VIE-based solvers, due to the extremely high contrast and high inhomogeneity of the scatterer.
	
	To demonstrate the superior accuracy of the proposed current-based VIE solver with PWL basis functions, we perform the numerical computation for the EM scattering for the "Duke" RHBM, where a high-permittivity dielectric pad is attached to the left side of the head of the model with dielectric properties of $\epsilon_r^{'} = 300$ and $\sigma = 0.25$ $\rm S / m$. The geometry and properties of the suspensions, as well as the incident field, which originates from a 16-rung high-pass birdcage coil with sinusoidal excitation tuned and matched at $f=300$ $\rm MHz$, are thoroughly explained in \cite{brink2016theoretical}.
	
\begin{figure}[ht!]
		\captionsetup{font=footnotesize}
		\centering
		\includegraphics[scale=0.78]{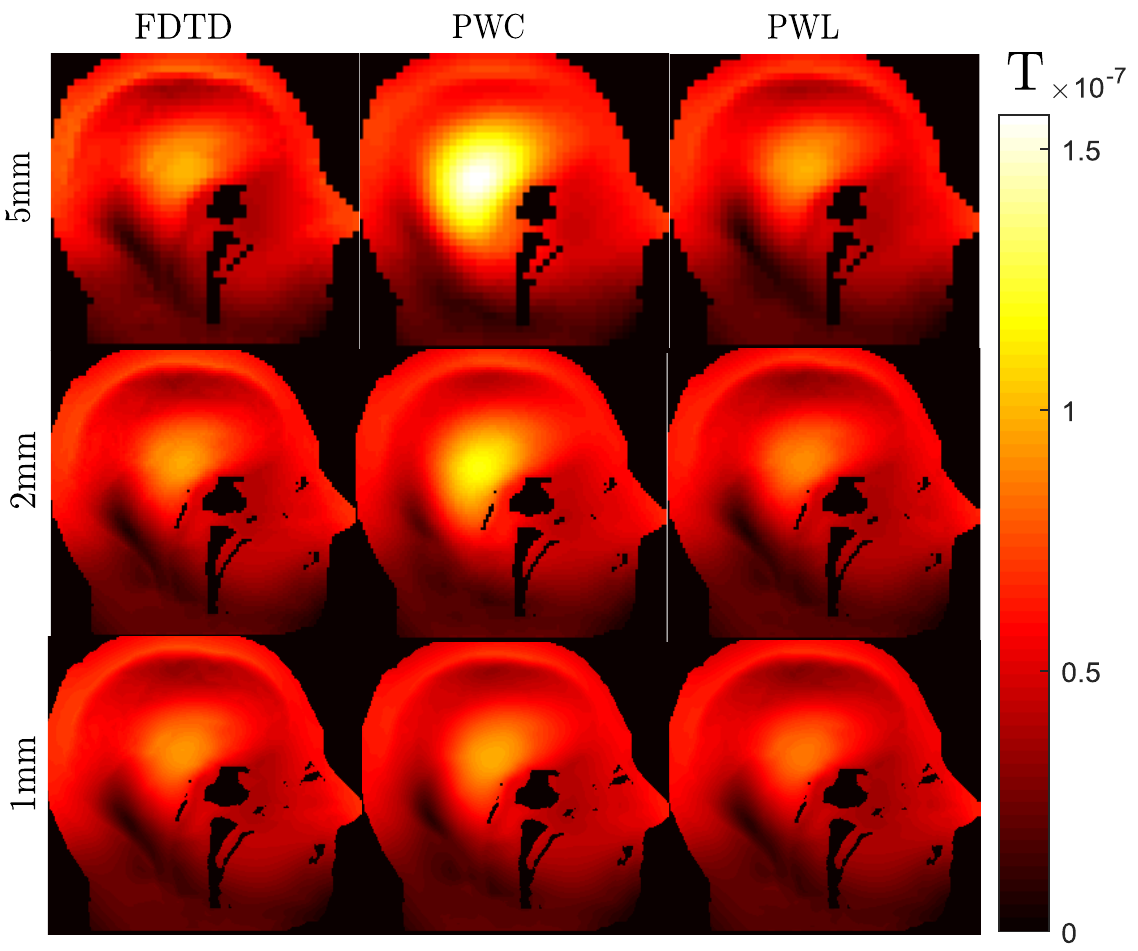}
		\caption{Sagittal views of $|b_1^+|$ in the "Duke" RHBM with an attached high-permittivity dielectric pad, when excited by a tuned birdcage coil at $7$ $\rm T$ MRI. From left to right, the magnetic field values are obtained with a commercial FDTD package, the PWC current-based solver, and the PWL current-based solver for resolutions  $h = 5$, $2$, and $1$ $\rm mm$, as viewed from top to bottom. Fields outside the body and at the pad are masked to enhance the visibility.}
		\label{jvie_vs_fdtd_b1}
	\end{figure}

	Subsequently, we calculate the EM field distributions over the head and the pad with the current-based VIE formulation with PWC and PWL basis functions, as we refine the resolution of the computational grid, using $h = 5$, $2$, and $1$ $\rm mm$. Table \ref{dofs_times_D} shows the DoFs and the memory footprint of each solver and the computational time that GMRES requires to converge with and without preconditioner. The matrix-vector product is performed in the GPU up to 2 mm for the PWC solver and up to 5 mm for PWL. The results are qualitatively compared to the EM fields obtained from a commercial software package which employs the FDTD method (xFDTD 7.2, Remcom Inc., State College, Pennsylvania, USA), as it is mentioned at \cite{brink2016theoretical}, and the $|b_1^+|$ at a sagittal cut for all $3$ resolutions and $3$ solvers is presented in Fig. \ref{jvie_vs_fdtd_b1}, masking the fields outside the body for improved visibility.
	
	\begingroup
	\setlength{\tabcolsep}{3pt}
	\begin{table}[ht!]
		\captionsetup{font=scriptsize}
		\caption{COMPUTATIONAL CONSIDERATIONS FOR EXPERIMENT D} \label{dofs_times_D} \centering
		\begin{tabular}{c|ccc}
			\hline\hline\\[-0.4em]
			$h(\rm mm)$  		        		& 5      	  & 2    		 & 1            \\[0.3em] \hline \\[-0.3em]
			$N_x\times N_y\times N_z$			& 39$\times$47$\times$45  & 97$\times$119$\times$112    & 194$\times$237$\times$224 \\[0.3em]
			DoFs PWC (million)					& 0.25   	  &  3.88         & 30.90       \\[0.3em]
			DoFs PWL (million)	        		& 0.99   	  & 15.51        & 123.59     \\[0.3em]\hline&&&\\[-0.4em]
			PWC Unknowns Mem. (GB)				& 0.029  	  & 0.462         &  3.683      \\[0.3em]
			PWL Unknowns Mem. (GB)	    		& 0.118  	  & 1.849    	  & 14.733     \\[0.3em]\hline&&&\\[-0.4em]
			PWC Time w/ prec					& 00:00:53    & 00:09:19     & 05:01:52     \\[0.3em]
			PWC Time w/o prec	        		& 00:01:14    & 00:19:55     & 09:10:35     \\[0.3em]
			PWL Time w/ prec	    			& 00:14:52    & 14:08:04     & 134:18:03    \\[0.3em]
			PWL Time w/o prec					& 01:12:19    & 84:07:55    & 488:01:24   \\[0.3em]\hline&&&\\[-0.4em]
			PWC  $W_{\rm m}(\rm nJ)$			& 2.394  	  & 1.932    	 & 1.689      \\[0.3em]
			PWL  $W_{\rm m}(\rm nJ)$	    	& 1.972  	  & 1.683    	 & 1.572    \\[0.3em]
			FDTD $W_{\rm m}(\rm nJ)$	    	& 2.356  	  & 2.013    	 & 1.893     \\[0.3em]\hline&&&\\[-0.4em]
			PWC $\textrm{error}_{\bm h_1}$ (\%) 			   	&52.33	 &22.94	   	&7.45		\\[0.3em]
			PWL $\textrm{error}_{\bm h_1}$ (\%)					&25.48	 & 7.10	   	&N/A			\\[0.3em]
			PWC $\textrm{error}_{\bm h_2}$ (\%)					&72.13	 &33.44	   	&13.42		\\[0.3em]
			PWL $\textrm{error}_{\bm h_2}$ (\%)					&36.21	 &10.34	    &N/A			\\[0.3em]
			PWC $\textrm{error}_{\bm h_3}$ (\%)					&83.84	 &49.76	    &24.86		\\[0.3em]
			PWL $\textrm{error}_{\bm h_3}$ (\%)					&35.49	 &15.51	    &N/A			\\[0.3em]
			\hline\hline
		\end{tabular}
	\end{table}
	\endgroup
	
	From that figure, it is clear that, even for this very demanding simulation, the PWL solver yields reliable and accurate EM fields for coarse resolutions, which are in good agreement with the FDTD-based solution, despite the inherent differences in the numerical modeling of the setting. On the contrary, it is necessary to refine the resolution up to $1$ $\rm mm$ for the PWC solver to converge to an accurate solution. Additionally, in Table \ref{dofs_times_D} we show the magnetic energy, defined as $W_m = 1/4\int_V\mu |\bm h|^2dV$, inside the head for all three resolutions and solvers and note the same decreasing trends with $h$-refinement for all solvers and with $p$-refinement for all resolutions. Additionally, the 2 mm PWL magnetic energy is lower than  the 1 mm PWC and the 5 mm PWL is approximately equal to the 2 mm PWC and under the same assumption that lower magnetic energies imply a lower errors, it could be argued that coarser PWL magnetic fields can provide more or equally accurate results with finer PWC magnetic fields that have more number of unknowns. Similarly to the previous section, we consider the finest 1 mm PWL magnetic fields as the reference solution and evaluate the following errors:
	\begin{subequations}
		\begin{gather}
		\textrm{error}_{\bm h_1} = |  W_{\rm m} - \hat{W}_{\rm m} |    / | \hat{W}_{\rm m} |, 		\\
		\textrm{error}_{\bm h_2} = || \mathbf w_{\rm m} - \hat{\mathbf w}_{\rm m} ||_1 / ||\hat{\mathbf w}_{\rm m} ||_1, 	\\  
		\textrm{error}_{\bm h_3} = || \mathbf h     - \hat{\mathbf h}	   ||_2  / ||\hat{\mathbf h} ||_2, 
		\end{gather}
	\end{subequations}
	where $\mathbf w_{\rm m} = 1/4\mu|\mathbf h|^2$ is the magnetic energy density. At Table \ref{dofs_times_D}, we show these errors and notice that PWL 2 mm is more accurate than PWC 1 mm for all these errors despite PWC having twice the unknowns of PWL. For the comparison between PWL 5 mm and PWC 2 mm, we see that PWL $\textrm{error}_{\bm h_3}$ is smaller than PWC $\textrm{error}_{\bm h_3}$ and PWL $\textrm{error}_{\bm h_1}$ and $\textrm{error}_{\bm h_2}$ are somewhat larger than the corresponding PWC, despite PWC having $\sim 3.75$ times the DoFs of PWL.
	
	\begin{figure}[ht!]
		\captionsetup{font=footnotesize}
		\centering
		\includegraphics[scale=0.3]{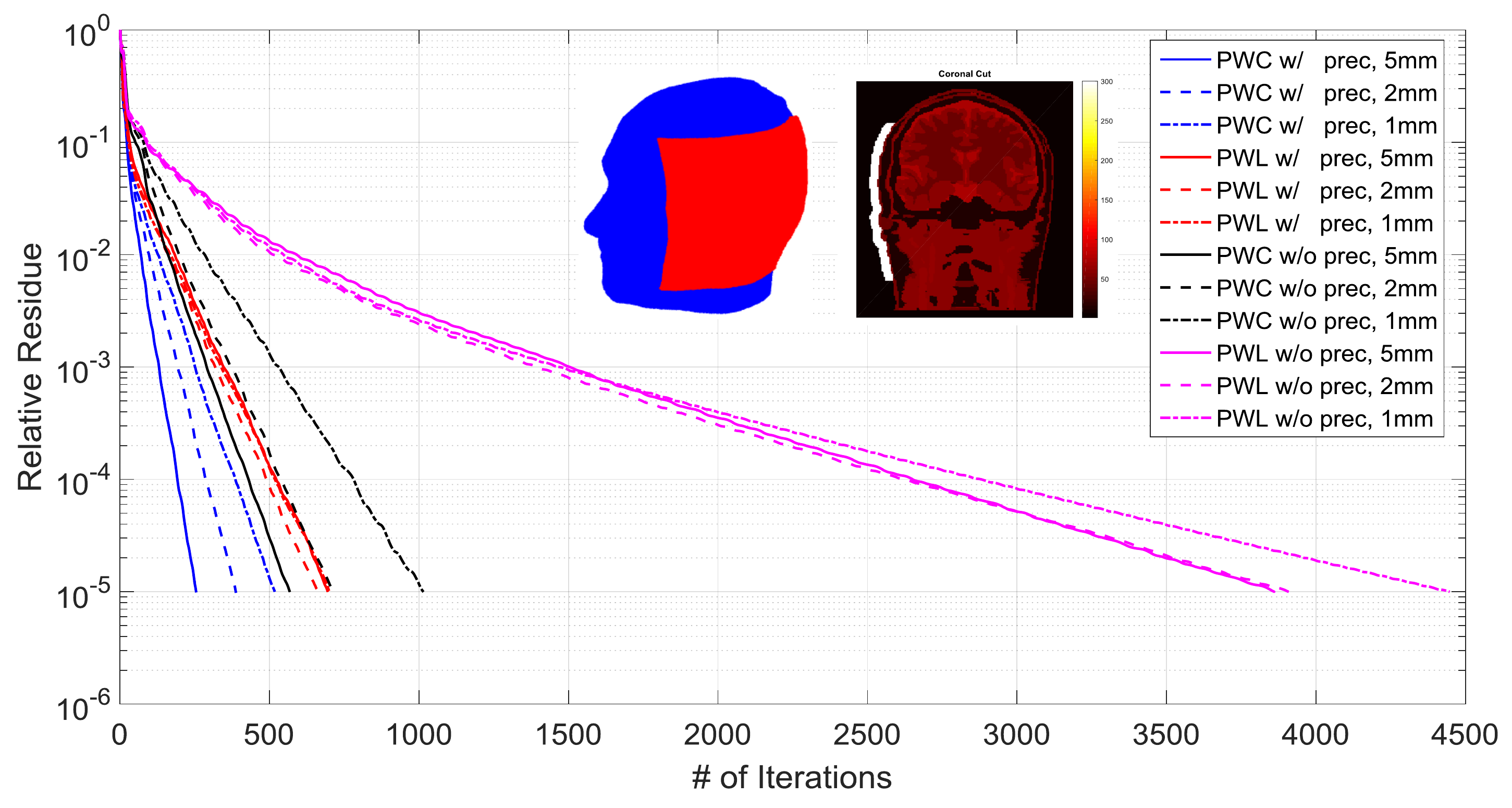}
		\caption{Convergence of GMRES iterative solver for the PWC current-based solver and the PWL current-based solver with and without preconditioner for different resolutions $h = 5$, $2$, and $1$ $\rm mm$, for the calculation of the EM scattering from the "Duke" RHBM with an attached high-permittivity dielectric pad irradiated by a tuned birdcage coil at $7$ $\rm T$ MRI. Also, the coronal view of the $\epsilon_r^{'}$ of the simulated geometry and the geometry of the attached pad are presented.}
		\label{convergence-jvie-pad}
	\end{figure}
	
	Finally, we execute the same numerical experiment for the PWC and PWL solvers utilizing the preconditioner mentioned in the previous section and present the convergence for the iterative solver for all the refinements in Fig. \ref{convergence-jvie-pad}. Here, the well-conditioned properties of the current-based VIE can be observed, even for this extremely challenging scenario, allowing both for the utilization of higher-order approximations and grid refinements. Similarly to the previous numerical experiment, the unpreconditioned PWL solver requires many iterations to converge. However, the effectiveness of the proposed preconditioner is again demonstrated, since when it is utilized the number of iterations drops from 3900 to 700, a value comparable to the iteration count of the PWC solver, allowing for this solver to be used in practical and demanding applications with manageable computational cost while simultaneously providing accurate and reliable results. Finally, it is worth noting that without the proposed preconditioner the iteration count of the PWL solver diverges at the finest resolution of $1$ $\rm mm$ which does not occur when the preconditioner is employed.
	
	\section{CONCLUSION} 
	A fast VIE solver based on the equivalent polarization/magnetization currents with PWL basis functions is derived for the accurate computation of the EM scattering from highly inhomogeneous and/or high contrast objects. The proposed solver has remarkably stable convergence properties and yields reliable EM fields for extremely challenging modeling scenarios and for coarse resolutions without necessarily refining the computational grid. Furthermore, by discretizing the VIE on uniform grids, the matrix-vector product can be performed fast with the help of FFT and, when combined with iterative solvers, large and complex problems can be solved accurately within reasonable time and computational resources. The proposed framework can be utilized in challenging applications such as the modeling of the interactions between EM fields and biological tissue, including the presence of shimming pads of very high electric permittivity.

	\section*{ACKNOWLEDGMENT}
	The authors would like to thank Wyger M. Brink for providing the FDTD data, the geometry of the attached pad, and for useful discussions. 
	
	\bibliographystyle{IEEEtran}
	\bibliography{IEEEabrv,FINAL_VERSION_ARXIV}

	\vskip -2\baselineskip plus -1fil
	\begin{IEEEbiographynophoto}
		{Ioannis P.~Georgakis} was born in Thessaloniki, Greece, in 1991. He received the Diploma degree in electrical engineering and computer science from the Aristotle University of Thessaloniki, Greece, in 2014 and the Ph.D. degree in computational and data science and engineering from the Skolkovo Institute of Science and Technology, Moscow, Russia, in 2019.
		\par
		Dr. Georgakis was a recipient of three scholarships of academic excellence by the Greek State Scholarship Foundation during his undergraduate studies (2010-2012). His research interests include computational electromagnetics, with an emphasis on volume integral equation methods for magnetic resonance imaging modeling.
	\end{IEEEbiographynophoto}
	\vskip -2\baselineskip plus -1fil
	
	\begin{IEEEbiographynophoto}
		{Ilias I.~Giannakopoulos}  was born in Athens, Greece, in 1993. He received the Diploma degree in electrical engineering and computer science from the Aristotle University of Thessaloniki, Greece, in 2016 and the Ph.D. degree in computational and data science and engineering from the Skolkovo Institute of Science and Technology, Moscow, Russia, in 2020.
		\par
		Dr. Giannakopoulos was a recipient of an honorary scholarship from the Greek State Scholarships Foundation (2014) for his excellence in his undergraduate studies for the academic year 2011-2012, and the IEEE AP-S Student Paper Competition Honorable Mention Award of the 2018 IEEE International Symposium on Antennas and Propagation. His research is focused on computational electromagnetics, with an emphasis on volume and surface integral equation methods, numerical linear algebra and magnetic resonance imaging.
	\end{IEEEbiographynophoto}
	\vskip -2\baselineskip plus -1fil
	
	\begin{IEEEbiographynophoto}
		{Mikhail S.~Litsarev} was born in Moscow Region, Russia, in 1983. He received the Diploma in condensed matter physics from Moscow Engineering Physics Institute in 2006 and the Ph.D. degree in theoretical physics from P.N. Lebedev Physical Institute of the Russian Academy of Sciences in 2010. 
		\par
		From 2011 to 2013, he was a Postdoctoral Research Associate at the Department of Physics and Astronomy, Uppsala University, Sweden, where he was a member of computational physics and condensed matter theory Group at the Materials Theory Division. From 2013 to 2018, he was a Postdoctoral Research Associate at the Skolkovo Institute of Science and Technology, Moscow, Russia. His research interests revolve around computational physics, numerical linear algebra and C++.
	\end{IEEEbiographynophoto}
	\vskip -2\baselineskip plus -1fil
	
	\begin{IEEEbiographynophoto}
		{Athanasios G.~Polimeridis} (SM'16) was born in Thessaloniki, Greece, in 1980. He received the Diploma and Ph.D. degrees in electrical engineering and computer science from the Aristotle University of Thessaloniki in 2003 and 2008, respectively.
		\par
		From 2008 to 2012, he was a Post-Doctoral Research Associate with the Laboratory of Electromagnetics and Acoustics, {\'E}cole Polytechnique F{\'e}d{\'e}rale de Lausanne, Switzerland. From 2012 to 2015, he was a Post-Doctoral Research Associate with the Massachusetts Institute of Technology, Cambridge, MA, USA, where he was a member of the Computational Prototyping Group, Research Laboratory of Electronics. From 2015 to 2018 he was an Assistant Professor with the Skolkovo Institute of Science and Technology, Moscow, Russia. He is currently the VP of Computation with Q Bio, CA, USA. His research interests include computational methods for problems in physics and engineering (classical electromagnetics, quantum and thermal electromagnetic interactions, and magnetic resonance imaging) with an emphasis on the development and implementation of integral-equation-based algorithms. 
		\par
		Dr. Polimeridis was a recipient of the Swiss National Science Foundation Fellowship for Advanced Researchers in 2012.
	\end{IEEEbiographynophoto}

\end{document}